\input amstex
\documentstyle{amsppt}

\topmatter
\title
An explicit formula for Hecke $L$-functions
\endtitle
\author
Xian-Jin Li
\endauthor
\abstract
In this paper an explicit formula is given for a sequence of 
numbers.  The positivity of this sequence of numbers implies 
that zeros in the critical strip of the Euler product of
 Hecke polynomials, which are associated with the space of 
 cusp forms of weight $k$ for Hecke congruence
 subgroups, lie on the critical line.
\endabstract
\address
Department of Mathematics, Brigham Young University,
Provo, Utah 84602
\endaddress
\email
xianjin\@math.byu.edu
\endemail
\thanks
Research supported by National Security Agency MDA904-03-1-0025
\endthanks
\endtopmatter
\document

\heading
1.  Introduction
\endheading

 Let $k$ and $N$ be positive integers with $k>2$, and let $\chi$ be
a Dirichlet character of modulus $N$ with $\chi(-1)=(-1)^k$
and with conductor $\frak f$.
We denote by $S_k(N, \chi)$ the space of all cusp forms
of weight $k$ and character $\chi$ for the Hecke congruence
subgroup $\Gamma_0(N)$ of level $N$.  That is,
$f$ belongs to $S_k(N, \chi)$ if and only if
$f$ is holomorphic in the upper half-plane, satisfies
$$f\left({az+b\over cz+d}\right)=\chi(d) (cz+d)^k f(z)$$
for all $\left(\smallmatrix a&b\\ c&d\endsmallmatrix\right)\in
\Gamma_0(N)$, satisfies the usual regularity conditions at the
cusps of $\Gamma_0(N)$, and vanishes at each cusp of
$\Gamma_0(N)$.

   The Hecke operators $T_n$ are defined by
$$(T_n f)(z)={1\over n}\sum_{ad=n}\chi (a)a^k
\sum_{0\leq b<d} f\left({az+b\over d}\right). \tag 1.1 $$
A function $f$ in $S_k(N, \chi)$ is called a
Hecke eigenform if
$$T_n f=\lambda(n)f$$
for all positive integers $n$ with $(n,N)=1$.
The Fricke involution $W$ is
defined by
$$(Wf)(z)=N^{-k/2}z^{-k}F(-1/Nz),$$
and the complex conjugation operator $K$ is defined by
$$(Kf)(z)=\bar f(-\bar z).$$
Set $\bar W=KW$.   An element $f$ in $S_k(N,\chi)$
is called a newform if it is an eigenfunction
of $\bar W$ and of all the Hecke operators $T_n$.
An element $f$ in $S_k(N,\chi)$ is called an oldform if there is an
element $f^\prime$ in $S_k(N^\prime, \chi_{N^\prime})$ such that
$f(z)=f^\prime(dz)$, where $N^\prime, d$ are positive integers
satisfying $\frak f|N^\prime$ and $dN^\prime|N$ and where
$\chi_{N^\prime}$ is the Dirichelt character of modulus $N^\prime$
induced by $\chi$.

   Let $f$ be a newform in $S_k(N, \chi)$ normalized so
that its first Fourier coefficient is $1$.  Then it
has the Fourier expansion
$$f(z)=\sum_{n=1}^\infty \lambda(n) e^{2\pi i nz}$$
with the Fourier coefficients equal to the eigenvalues of
Hecke operators.  Since $f$ is an eigenfunction of the involution
 $\bar W$, we can assume that
$$\bar Wf=\eta f$$
for constant $\eta$.  Let
$$L_f(s)=\sum_{n=1}^\infty {\lambda(n)\over n^s}$$
for $\Re s>{k+1\over 2}$.
It has the Euler product
$$L_f(s)=\prod_p (1-\lambda(p)p^{-s}+\chi(p)p^{k-1-2s})^{-1}.
\tag 1.2$$
If we denote
$$\xi_f(s)=\left(\sqrt N\over 2\pi\right)^s
\Gamma({k-1\over 2}+s)L_f({k-1\over 2}+s), \tag 1.3$$
then $\xi_f(s)$ is an entire function and satisfies the
functional identity
$$\xi_f(s)=i^k\bar \eta\bar\xi_f(1-\bar s).$$
 When $\chi$ is primitive, we have
 $\eta=\tau(\bar\chi)\lambda(N)N^{-k/2}$ with
$\tau(\chi)$ being the Gauss sum for $\chi$; see Iwaniec [6].

 Let
$$L_N(s)=\prod_{p\nmid N}\det|1-T(p)p^{-s}+\chi(p)p^{k-1-2s}I|^{-1}
\tag 1.4$$
where $I$ is the identity map acting on the space $S_k(N, \chi)$,
and let
$$\xi_N(s)=N^{gs/2}(2\pi)^{-gs}\Gamma^g({k-1\over 2}+s)
L_N({k-1\over 2}+s) \tag 1.5$$
where $g$ is the dimension of the space $S_k(N,\chi)$.
We will show that $\xi_N(s)$ is an entire function, and
that zeros of $\xi_N(s)$ in the strip $0<\Re s<1$ appear in
pairs $\rho$ and $1-\bar\rho$.

   Let
$$\tau_N(n)=\sum_\rho [1-(1-{1\over\rho})^{-n}] \tag 1.6 $$
for $n=1,2,\cdots$, where the sum on $\rho$ runs over all
zeros of $\xi_N(s)$ taken in the order given by
$|\Im \rho|<T$ for $T\to\infty$ with a zero of multiplicity
$\ell$ appearing $\ell$ times in the list.  If $\rho=0$ is
a zero of $\xi_N(s)$, then $(1-1/\rho)^{-n}$ in (1.6) is
interpreted to be $0$.

In [2], Bombieri and Lagarias generalized a criterion of
the author for the Riemann hypothesis [8] and obtained
the following useful theorem.

\proclaim{Theorem 1.1} (Bombieri-Lagarias [2])
 Let $\Cal R$ be a set of complex
numbers $\rho$, whose elements have positive integral multiplicities
assigned to them, such that $1\not\in\Cal R$ and
$$\sum_\rho {1+|\Re \rho|\over (1+|\rho|)^2}<\infty.$$
Then the following conditions are equivalent:
\roster
\item $\Re\rho\leq {1\over 2}$ for every $\rho$ in $\Cal R$;
\item $\sum_\rho \Re[1-(1-{1\over\rho})^{-n}]\geq 0$
   for $n=1,2,\cdots$.
\endroster
\endproclaim

   As a corollary of Theorem 1.1 we obtain a criterion
for the location of all nontrivial zeros of Hecke $L$-functions
associated with all Hecke eigenforms which form an orthonormal basis
in $S_k(N,\chi)$.

\proclaim{Corollary 1.2}  All zeros of $\xi_N(s)$ in the
 strip $0<\Re s<1$ lie on the critical line
$\Re s=1/2$ if, and only if, $\tau_N(n)\geq 0$
for all positive integers $n$.
\endproclaim

  In order to state the main theorem, we need an
explicit formula for the trace tr$T(p^\ell)$ of Hecke
operators $T(p^\ell)$ acting on the space $S_k(N,\chi)$
for all primes $p\nmid N$ and for $\ell=1,2,\cdots$.
This formula is given by the Eichler-Selberg
trace formula obtained in Oesterl\'e [12] (cf. Serre [13]).
We denote by $\varphi$ the Euler
$\varphi$-function.  Let
$$\psi(N)=N\prod_{p|N}  (1+1/p),$$
and let $\chi(\sqrt n)=0$ if $n$ is not the square of an integer.

  \proclaim{Theorem 1.3} (Th\'eor\`eme $3^\prime$, [12]) For every
  positive integer $n$, the trace tr$(T(n))$ of the Hecke operator $T(n)$
  acting on the space $S_k(N,\chi)$ is given by
  $$\aligned \text{tr}(T(n))&=n^{{k\over 2}-1}\chi(\sqrt n)
  {k-1\over 12}\psi(N) \\
  &-\sum_{t\in\Bbb Z, t^2<4n}{\rho^{k-1}-\bar\rho^{k-1}
  \over \rho-\bar\rho}\sum_{\underset{{t^2-4n\over m^2}
  \equiv 0, 1 (\mod 4)}\to {m\in \Bbb Z^+}}{h((t^2-4n)/m^2)
  \over w((t^2-4n)/m^2)}\mu(t, n, m)\\
  &-{1\over 2}\sum_{0<d|n}\min(d^{k-1}, (n/d)^{k-1})
  \sum_{\underset{(c, {N\over c})|
  ({N\over\frak f}, {n\over d}-d)}\to {c|N}}\varphi((c, N/c))\chi(y)
   \endaligned $$
  where $\rho, \bar\rho$ are the roots of $x^2-tx+n=0$,
  where the integer $y$ modulo $N/(c, N/c)$ is defined by
  $y\equiv d (\mod c), y\equiv n/d (\mod N/c)$, where
  $$\mu(t, n, m)={\psi(N)\over \psi(N/\gcd(N,m))}
  \sum_{\underset{x^2-tx+n\equiv 0 (\mod N(N, m))}\to
  {x (\mod N)}}\chi(x),$$
  and where $h(f)$ and $w(f)$ are respectively the class number
  and the number of units in the ring of integers of the
  imaginary quadratic field of discriminant $f<0$.
   \endproclaim

  Let $\gamma=0.5772\cdots$ be Euler's constant, let
$$\Lambda(m)=\cases \ln p, &\text{if $m$ is a
positive power of a prime $p$;}\\
0, &\text{otherwise,} \endcases $$
and let $d(m)$ be the number of positive divisors
of $m$.  For $m|N$ and $\frak f|m$, we denote by $\nu_m$
the dimension of the subspace generated by all newforms
in $S_k(m, \chi_m)$, where $\chi_m$ is the
Dirichlet character of modulus $m$ induced by the Dirichlet
character $\chi$ of modulus $N$.  In this paper we obtain the
following explicit formula for the $\tau_N(n)$'s.
The weight $2$ case of this formula was obtained in Li [9].

 \proclaim{Theorem 1.4}  Let $\tau_N(n)$ be given in
 (1.6).  Then we have
$$\aligned \tau_N(n)=& {n\over 2}
\ln\left(N^{\nu_N}\prod_{\frak f|m, m|N} m^{\nu_md(N/m)}\right)
-ng(\ln 2\pi +\gamma+{2\over k+1})\\
&-\sum_{l=1}^n \binom nl {(-1)^{l-1}\over (l-1)!}
\sum_{\underset{(m, N)=1}\to {m=1}}^\infty
{\Lambda(m)\over m^{k+1\over 2}}B(m)(\ln m)^{l-1}\\
&+ng\sum_{l=1}^\infty {k+1\over l(2l+k+1)}
+g\sum_{m=2}^n\binom nm \sum_{l=1}^\infty
{(-1)^m\over (l+{k-1\over 2})^m}. \endaligned $$
for all positive integers $n$, where
$B(p^\ell)= \text{tr}(T(p^\ell))-\chi(p)p^{k-1}
\text{tr}(T(p^{\ell-2}))$ for $p\nmid N$.
\endproclaim

   This paper is organized as follows:   An arithmetic formula
is obtained in section 2 for a sequence of numbers, whose
positivity implies that nontrivial zeros of Hecke $L$-functions
associated with newforms lie on the critical line.
 This formula will be used in the proof of
 Theorem 1.4.   In section 3 we give some preliminary
 results for the proof of Theorem 1.4.  Finally,
 Theorem 1.4 is proved in section 4.

   This research started while the author attended the Workshop on
Zeta-Functions and Associated Riemann Hypotheses, New York
University, Manhattan, May 29 - June 1, 2002.  He wants to thank
the American Institute of Mathematics, Brian Conrey, and Peter
Sarnak for the support.

\heading
2.  Explicit formulas for Hecke $L$-functions
\endheading

   Let $f$ be a normalized newform in $S_k(N,\chi)$, and let
$\xi_f(s)$ be given in (1.3).  Put
$$\tau_f(n)=\sum_\rho [1-(1-{1\over\rho})^{-n}] \tag 2.1 $$
for $n=1,2,\cdots$, where the sum is over all the zeros of $\xi_f(s)$
in the order given by $|\Im\rho|<T$ for $T\to\infty$
with a zero of multiplicity
$\ell$ appearing $\ell$ times in the list.

   Assume that $f$ is a normalized newform in $S_k(N,\chi)$.
For each prime number $p$, let $\alpha_p$ and $\beta_p$
be the two roots of $T^2-\lambda(p)T+\chi(p)p^{k-1}$ where $\lambda(p)$
is given in (1.2).  Put
$$b_f(p^m)=\cases\lambda(p)^m, &\text{ if $p|N$;}\\
\alpha_p^m+ \beta_p^m, &\text{ if $(p,N)=1$.}\endcases
\tag 2.2$$

   The flowing arithmetic formula for $\tau_f(n)$ generalizes
an arithmetic formula of Bombieri and Lagarias [2] for the
Riemann zeta function.

  \proclaim{Theorem 2.1}  Assume that $f$ is a normalized newform
in $S_k(N, \chi)$.  If $\tau_f(n)$ is given in (2.1), then we have
$$ \aligned \tau_f(n)=& n\left(\ln {\sqrt N\over 2\pi}-\gamma\right)
-\sum_{j=1}^n \binom nj {(-1)^{j-1}\over (j-1)!}\sum_{l=1}^\infty
{\Lambda(l)\over l^{k+1\over 2}} b_f(k)(\ln l)^{j-1}\\
&+n\left(-{2\over k+1}+\sum_{l=1}^\infty {k+1\over l(2l+k+1)}\right)
+\sum_{j=2}^n\binom nj(-1)^j\sum_{l=1}^\infty
{1\over (l+{k-1\over 2})^j}. \endaligned $$
for $n=1,2,\cdots$. \endproclaim

   \proclaim{Lemma 2.2} (see [10])
Let $F(x)$ be a function defined on the real line $\Bbb R$ such that
$$2F(x)=F(x+0)+F(x-0)$$
for all $x\in\Bbb R$, such that $F(x)\exp((\epsilon+1/2)|x|)$
is integrable and of bounded variation on $\Bbb R$ for a
constant $\epsilon >0$, and such that
$(F(x)-F(0))/x$ is of bounded variation on $\Bbb R$.   Then
$$\aligned \sum_\rho\Phi(\rho)
= &2F(0)\ln{\sqrt N\over 2\pi}-\sum_{n=1}^\infty
{\Lambda(n)\over n^{k/2}}[b_f(n)F(\ln n)+\bar b_f(n)F(-\ln n)]\\
&-\int_0^\infty \left((F(x)+F(-x)){e^{-kx/2}
\over 1-e^{-x}}-2F(0) {e^{-x}\over x}\right)dx,\endaligned $$
where the sum on $\rho$ runs over all
zeros of $\xi_f(s)$ in the order given by $|\Im\rho|<T$ for
$T\to\infty$, and
$$\Phi(s)=\int_{-\infty}^\infty F(x) e^{(s-1/2)x}dx. $$
\endproclaim

\proclaim{Lemma 2.3} ([4] [5] [11])  Let $f$ be a newform in
$S_k(N, \chi)$.  Then there an absolute effective constant $c>0$
such that $L_f({k-1\over 2}+s)$ has no zeros in the region
$$\{s=\sigma+it:\, \sigma\geq 1-{c\over \ln(N+1+|t|)}\},$$
where $L_f$ is given in (1.2).
\endproclaim

 \proclaim{Lemma 2.4}  (Lemma 2 of [2])  For $n=1,2,\cdots$, let
$$F_n(x)=\cases e^{x/2}\sum_{j=1}^n \binom nj {x^{j-1}\over (j-1)!},
&\text{if $-\infty<x<0$};\\
n/2, &\text{if $x=0$};\\
0, &\text{if $0<x$}. \endcases $$
Then
$$\Phi_n(s)=1-\left(1-{1\over s}\right)^n$$
where $\Phi_n$ is related to $F_n$ by the relation
$$\Phi_n(s)=\int_{-\infty}^\infty F_n(x) e^{(s-1/2)x}dx. $$
\endproclaim

  \demo{Proof of Theorem 2.1}
 Since $\xi_f(s)$ is an entire function of order one
  and satisfies the functional identity
  $\xi_f(s)=i^k\bar\eta\bar\xi_f(1-\bar s)$,  we have
$$\xi_f(s)=i^k\bar\eta\bar\xi_f(1)\prod_\rho (1- s/\rho)$$
where the product is over all the zeros of $\xi_f(s)$ in
the order given by $|\Im\rho|<T$ for $T\to\infty$.
If $\varphi_f(z)=\xi_f(1/(1-z))$, then
$${\varphi_f^\prime(z)\over\varphi_f(z)}=\sum_{n=0}^\infty
\tau_f(n+1)z^n \tag 2.3$$
where the coefficients $\tau_f(n)$ are given in (2.1).

    For a sufficiently large positive number $X$
that is not an integer, let
$$F_{n,X}(x)=\cases F_n(x), &\text{if $-\ln X<x<\infty$};\\
{1\over 2}F_n(-\ln X), &\text{if $x=-\ln X$}; \\
0, &\text{if $-\infty<x<-\ln X$} \endcases $$
where $F_n(x)$ is given in Lemma 2.4.  Then $F_{n,X}(x)$
satisfies all conditions of Lemma 2.2.  Let
$$\Phi_{n, X}(s)=\int_{-\infty}^\infty F_{n, X}(x)
e^{(s-1/2)x}dx.$$
By Lemma 2.2, we obtain that
$$\aligned \sum_\rho \Phi_{n,X}(\rho)
= &2F_{n,X}(0)\ln {\sqrt N\over 2\pi}
-\sum_{l=1}^\infty {\Lambda(l)\over l^{k/2}} [b_f(l)
F_{n,X}(\ln l)+\bar b_f(l)F_{n,X}(-\ln l)]\\
&-\int_0^\infty \left({F_{n,X}(x)+F_{n,X}(-x)\over 1-e^{-x}}
e^{-kx/2}-2F_{n,X}(0)
{e^{-x}\over x}\right)dx,\endaligned $$
where the sum on $\rho$ runs over all zeros of $\xi_f(s)$ in the
order given by $|\Im\rho|<T$ for $T\to\infty$.
It follows that
$$\aligned \lim_{X\to\infty}&\sum_\rho\Phi_{n,X}(\rho)\\
&= n\left(\ln {\sqrt N\over 2\pi}-\gamma\right)
-\sum_{j=1}^n \binom nj {(-1)^{j-1}\over (j-1)!}\sum_{l=1}^\infty
{\Lambda(l)\over l^{k+1\over 2}} \bar b_f(k)(\ln l)^{j-1}\\
&+n\left(-{2\over k+1}+\sum_{l=1}^\infty {k+1\over l(2l+k+1)}\right)
+\sum_{j=2}^n\binom nj(-1)^j\sum_{l=1}^\infty
{1\over (l+{k-1\over 2})^j}. \endaligned $$
 We have
$$\aligned \Phi_n(s)-\Phi_{n,X}(s)&=X^{-s}\sum_{j=1}^n\binom nj
(-1)^{j-1}\sum_{m=0}^{j-1}{(\ln X)^{j-m-1}\over (j-m-1)!}
 s^{-m-1}\\
 &={X^{-s}\over s}\sum_{j=1}^n\binom nj {(-\ln X)^{j-1}\over (j-1)!}
 +O\left({(\ln X)^{n-2}\over |s|^2}
 X^{-\Re s}\right) \endaligned  \tag 2.4$$
for $\Re s>0$.

   Let $\rho$ be any zero of $\xi_f(s)$.
 By Lemma 2.3, we have
 $${c\over \ln (N+1+|\rho|) }
 \leq \Re\rho\leq 1-{c\over \ln (N+1+|\rho|)}  $$
 for a positive constant $c$.  An argument similar to that
 made in the proof of (3.9) of [2] shows that
 $$\sum_\rho {X^{-\Re\rho}\over |\rho|^2}
 \ll e^{-c^\prime\sqrt{\ln X}} \tag 2.5$$
 for a positive constant $c^\prime$.

    Since
     $$\aligned \sum_\rho {X^{-\rho}\over \rho}&=\sum_\rho
   {X^{-(1-\bar\rho)}\over 1-\bar\rho} \\
   &=-{1\over X}\sum_\rho {X^{\bar\rho}\over\bar\rho}
   +O\left(\sum_\rho {X^{-(1-\Re\rho)}\over |\rho|^2}\right)\\
   &=-{1\over X}\sum_\rho {X^{\bar\rho}\over\bar\rho}
 +O\left(e^{-c^\prime\sqrt{\ln X}}\right),
   \endaligned $$
and since
$$\lim_{X\to\infty}{(\ln X)^{j-1}\over X}\sum_\rho
{X^{\bar\rho}\over\bar\rho}=0 $$
for $j=1,2,\cdots, n$ by Theorem 4.2 and Theorem 5.2 of [11],
we have
$$\lim_{X\to\infty}(\ln X)^{j-1}\sum_\rho {X^{-\rho}\over \rho}=0 \tag 2.6$$
 for $j=1,2,\cdots, n$.  It follows from (2.4), (2.5) and
 (2.6) that
  $$\lim_{X\to\infty}\sum_\rho\Phi_{n,X}(\rho)
  =\sum_\rho\Phi_n(\rho).$$
 Hence, we have
$$\aligned \sum_\rho\Phi_n(\bar\rho)
= & n\left(\ln {\sqrt N\over 2\pi}-\gamma\right)
-\sum_{j=1}^n \binom nj {(-1)^{j-1}\over (j-1)!}\sum_{l=1}^\infty
{\Lambda(l)\over l^{k+1\over 2}} b_f(k)(\ln l)^{j-1}\\
&+n\left(-{2\over k+1}+\sum_{l=1}^\infty {k+1\over l(2l+k+1)}\right)
+\sum_{j=2}^n\binom nj(-1)^j\sum_{l=1}^\infty
{1\over (l+{k-1\over 2})^j}. \endaligned $$
Since zeros of $\xi_f(s)$ appear in pairs $\rho$ and $1-\bar \rho$,
we have
$$ \aligned \sum_\rho\Phi_n(\bar\rho)
=&\sum_\rho [1-(1-{1\over\bar\rho})^n] \\
&=\sum_\rho [1-(1-{1\over 1-\rho})^n\\
&=\sum_\rho [1-(1-1/\rho)^{-n}]=\tau_f(n)\endaligned $$
for $n=1,2,\cdots$, where the sum is over all the zeros of $\xi_f(s)$
in the order given by $|\Im\rho|<T$ for $T\to\infty$
with a zero of multiplicity
$\ell$ appearing $\ell$ times in the list.

This completes the proof of the theorem.
\qed\enddemo

\heading
3.  Preliminary results
\endheading

  A fundamental result of Hecke asserts that a basis
$\{f_1, f_2, \cdots, f_g\}$ in $S_k(N, \chi)$ exists which consists
of eigenfunctions of all the Hecke operators $T(n)$ with
$(n, N)=1$;  see  Theorem 6.21 in Iwaniec [6].   We can assume
that each $f_j$ is either a normalized newform in $S_k(N, \chi)$ or
coming from a normalized newform in a lower level.
For $j=1,\cdots, g$,
we choose $g_j=f_j$ if $f_j$ is a normalized newform in $S_k(N, \chi)$,
and $g_j=f^\prime_j$ if $f_j$ is an oldform in $S_k(N, \chi)$
and if $f^\prime_j$ is a normalized newform in
$S_k(N^\prime, \chi_{N^\prime})$ for
some divisor $N^\prime$ of $N$ with $\frak f|N^\prime$ such that
$f_j(z)=f^\prime_j(dz)$ for some positive integer $d|N/N^\prime$,
where $\chi_{N^\prime}$ is the Dirichlet character of modulus
$N^\prime$ induced by the Dirichlet character $\chi$ of modulus $N$.

   Let
$$\xi_H(s)=\prod_{j=1}^g \xi_{g_j}(s) \tag 3.1$$
where $\xi_{g_j}(s)$ is defined as in (1.3).
Since $\xi_{g_j}(s)$ is an entire function and satisfies
the functional identity $\xi_{g_j}(s)=w_j \bar\xi_{g_j}(1-\bar s)$
for a constant $w_j$, the function $\xi_H(s)$ is entire
and satisfies the functional identity
$$\xi_H(s)=w \bar\xi_H(1-\bar s) \tag 3.2$$
for a constant $w$.  Put
$$\tau_H(n)=\sum_{j=1}^g\tau_{g_j}(n),\tag 3.3$$
where $\tau_{g_j}(n)$ is defined similarly as in (2.1).
If  $\varphi_H(z)=\xi_H(1/(1-z))$, then we have
$${\varphi_H^\prime(z)\over\varphi_H(z)}=\sum_{n=0}^\infty
\tau_H(n+1)z^n \tag 3.4$$
by (2.3).

  \proclaim{Lemma 3.1}  For all positive integers $n$ with
  $(n, N)=1$, we have
    $$(T(n)f_j)(z)=\lambda_{g_j}(n)f_j(z)$$
    for $j=1,2,\cdots, g$, where $\lambda_{g_j}(n)$ is
  the eigenvalue of $T(n)$ acting on $g_j(z)$. \endproclaim

\demo{Proof}  If $f_j$ is a newform, the stated identity is
trivially true.

   Next, we assume that $f_j$ is an oldform.
Let $g_j$ be a normalized newform in $S_k(N^\prime,\chi_{N^\prime})$
for some divisor $N^\prime$ of $N$ with $\frak f|N^\prime$
such that $f_j(z)=g_j(dz)$
for some positive integer $d|N/N^\prime$.  Since $(n,N)=1$,
by (1.1) we have
$$(T(n) f)(z)={1\over n}\sum_{ad=n}\chi(a)a^k
\sum_{0\leq b<d} f\left({az+b\over d}\right)$$
for any function $f$ in $S_k(N, \chi)$.  Thus, we have
$$(T(n) f_j)(z)={1\over n}\sum_{\alpha\delta=n}
\chi(\alpha)\alpha^k\sum_{0\leq \beta<\delta} g_j
\left({\alpha dz+d\beta\over \delta}\right). $$
Since $(n, N)=1$, $d|N$ and $\alpha\delta=n$, we have
$(\delta, d)=1$.  Let $r_\beta$ be remainder of
$d\beta$ modulo $\delta$.  Then $\{r_\beta: 0\leq\beta<\delta\}
=\{0, 1, \cdots, \delta-1\}$.  Since
$g_j\in S_k(N^\prime,\chi_{N^\prime})$, we have
$$g_j\left({\alpha dz+d\beta\over \delta}\right)
=g_j \left({\alpha dz+r_\beta\over \delta}\right).$$
Note that $\chi(t)=\chi_{N^\prime}(t)$ for any positive integer
$t$ with $(t, N)=1$.  It follows that
$$\aligned (T(n) f_j)(z)&={1\over n}\sum_{\alpha\delta=n}
\chi_{N^\prime}(\alpha)\alpha^k\sum_{0\leq \beta<\delta} g_j
\left({\alpha dz+\beta\over \delta}\right)\\
&=(T(n) g_j)(w)=\lambda_{g_j}(n)g_j(w)
=\lambda_{g_j}(n)f_j(z) \endaligned  $$
where $w=dz$.

  This completes the proof of the lemma.
\qed\enddemo

  \proclaim{Lemma 3.2}  Let $\xi_N(s)$ be given in (1.5).
  Then $\xi_N(s)$ is an entire function, and its zeros
 in the strip $0<\Re s<1$ appear in pairs $\rho$ and $1-\bar \rho$.
 \endproclaim

  \demo{Proof}    For $j=1,2, \cdots, g$, let $g_j$ be a normalized
  newform in $S_k(N_j,\chi_j)$ where $\chi_j$ is the Dirichlet
character of modulus $N_j$ induced by the Dirichlet character
$\chi$ of modulus $N$.  Then we can write
  $$L_{g_j}(s)=\prod_p (1-\lambda_{g_j}(p)
  p^{-s}+\chi_j(p)p^{k-1-2s})^{-1}. \tag 3.5$$
Since $f_1, f_2, \cdots, f_g$ are eigenfunctions of all the
Hecke operators $T(n)$ with $(n, N)=1$, by Lemma 3.1 we have
$$\det|1-T(p)p^{-s}+\chi(p)p^{k-1-2s}I|=\prod_{j=1}^g
(1-\lambda_{g_j}(p)p^{-s}+\chi(p)p^{k-1-2s})$$
for any prime $p\nmid N$.  Let $L_N(s)$ be given in (1.4).
Then we have
$$L_N(s) =\prod_{j=1}^g\prod_{p\nmid N}
(1-\lambda_{g_j}(p)p^{-s}+\chi(p)p^{k-1-2s})^{-1}.$$
Since $\chi(p)=\chi_j(p)$ for $j=1,\cdots, g$ when $p\nmid N$,
we have
$$L_N(s)=\prod_{j=1}^g(L_{g_j}(s)\prod_{p|N_j}(1-\lambda_{g_j}(p)p^{-s})
\prod_{p\nmid N_j,p|N}(1-\lambda_{g_j}(p)p^{-s}+\chi_j(p)p^{k-1-2s})). $$
Since
$$\xi_{g_j}(s)=N_j^{s/2}(2\pi)^{-s}
\Gamma({k-1\over 2}+s)L_{g_j}({k-1\over 2}+s),$$
we have
$$\xi_H(s)=A^{s/2}N^{gs/2}(2\pi)^{-gs}\Gamma^g
\left({k-1\over 2}+s\right)\prod_{j=1}^gL_{g_j}
\left({k-1\over 2}+s\right)$$
where $A=N^{-g}\prod_{j=1}^g N_j$.  It follows that
$$\aligned \xi_N(s) =&\xi_H(s) A^{-s/2} \\
&\times \prod_{j=1}^g(\prod_{p|N_j}(1-\lambda_{g_j}(p)p^{{1-k\over 2}-s})
\prod_{p\nmid N_j,p|N}(1-\lambda_{g_j}(p)p^{{1-k\over 2}-s}+\chi_j(p)p^{-2s})).
\endaligned \tag 3.6$$
This implies that $\xi_N(s)$ is an entire function.

 Assume $p|N_j$.  By Theorem 3 of Li [7], if $\chi_j$ is not a character
mod $N_j/p$ then $|\lambda_{g_j}(p)|=p^{k-1\over 2}$, and
if $\chi_j$ is a character mod $N_j/p$ then $\lambda_{g_j}(p)=0$
when $p^2|N_j$ and $|\lambda_{g_j}(p)|=p^{k-2\over 2}$ when $p^2\nmid N_j$.
It follows that all zeros of the function
$$ \prod_{j=1}^g\prod_{p|N_j}\left(1-\lambda_{g_j}(p)p^{{1-k\over 2}-s}\right)$$
lie on the lines $\Re s=0, -1/2$.

 Assume that $p\nmid N_j$ and $p|N$.  Since
the two roots of the polynomial
$1-\lambda_{g_j}(p)p^{1-k\over 2}z+\chi_j(p)z^2$ for $p\nmid N_j$
are of absolute value one by
the Ramanujan conjecture which was proved in Th\'eor\`eme 8.2 of
Deligne [3], all zeros of the function
$$ \prod_{j=1}^g \prod_{p\nmid N_j,p|N}\left(1-\lambda_{g_j}(p)
p^{{1-k\over 2}-s}+\chi_j(p)p^{-2s}\right)$$
lie on the line $\Re s=0$.
Therefore, it follows from (3.2) and (3.6) that
zeros of $\xi_N(s)$ in the critical strip $0<\Re s<1$
 appear in pairs $\rho, 1-\bar\rho$.

  This completes the proof of the lemma.
\qed\enddemo

\demo{Proof of Corollary 1.2}  Let $\tau_N(n)$ be defined by
(1.6) for all positive integers $n$.
Since $\xi_H(s)$ is an entire function of order one,
by (3.6) $\xi_N(s)$ is an entire function of order one.
This implies that
$$\sum_\rho {1+|\Re \rho|\over (1+|\rho|)^2}<\infty,$$
where the sum is over all zeros $\rho$ of $\xi_N(s)$.
Thus, conditions of Theorem 1.1 are satisfied.
Since by (3.6) all zeros of $\xi_N(s)$ outside the
strip $0<\Re s<1$ lie on the lines $\Re s=0, -1/2$,
 Theorem 1.1 implies that all zeros of $\xi_N(s)$ in the
 critical strip $0<\Re s<1$ satisfy
$\Re s\leq 1/2$ if, and only if, $\tau_N(n)\geq 0$
for all positive integers $n$.  By Lemma 3.2, all zeros of
of $\xi_N(s)$ in the  critical strip $0<\Re s<1$
appear in pairs $\rho$ and $1-\bar\rho$.
Thus, $\tau_N(n)\geq 0$
for all positive integers $n$ if, and only if,
$\Re \rho\leq 1/2$ and $\Re (1-\bar\rho)\leq 1/2$
for all zeros $\rho$ of $\xi_N(s)$ in the
critical strip $0<\Re s<1$.  That is,
all zeros of $\xi_N(s)$ in the
 strip $0<\Re s<1$ lie on the critical line
$\Re s=1/2$ if, and only if, $\tau_N(n)\geq 0$
for all positive integers $n$.

 This completes the proof of the corollary.
\qed\enddemo

\proclaim{Lemma 3.3}  Let $p$ be a prime, and let
$\alpha$ be a complex number.  Then we have
$$1-\alpha p^{-s}=c_p s^{\epsilon_p}\prod_\rho (1-s/\rho)$$
where the product on $\rho$ is over all nonzero zeros
of $1-\alpha p^{-s}$
taken in the order given by $|\rho|<T$ for $T\to\infty$ and where
$c_p=1-\alpha, \epsilon_p =0$ if $\alpha\neq 1$ and
$c_p=\ln p, \epsilon_p=1$ if $\alpha=1$.
\endproclaim

 \demo{Proof}  Since $1-\alpha p^{-s}$ is an entire function of
 order one, by Hadamard's factorization theorem there
 is a constant $a$ such that
 $$1-\alpha p^{-s}=c_pe^{as}s^{\epsilon_p}
 \prod_\rho (1-s/\rho)e^{s/\rho} \tag 3.7$$
 where the product is over all nonzero zeros of $1-\alpha p^{-s}$.
 Let $\alpha=e^{it}$ with $0\leq \Re t<2\pi$.
 Then the zeros of $1-\alpha p^{-s}$
are $i(t+2k\pi)/\ln p$, $k=0, \pm 1, \pm 2, \cdots$.
Since
$$\sum_{k=1}^\infty \left({-i\ln p\over t+2k\pi}
+{-i\ln p\over t-2k\pi}\right)
=2it\ln p\sum_{k=1}^\infty {1\over (2k\pi)^2-t^2} $$
is absolutely convergent,  by using (3.7) we can write
$$1-\alpha p^{-s}=c_pe^{hs}s^{\epsilon_p}
 \prod_\rho (1-s/\rho) \tag 3.8 $$
 for a constant $h$, where the product runs over all
 nonzero zeros $\rho$ of $1-\alpha p^{-s}$
taken in the order given by $|\rho|<T$ for $T\to\infty$.
By taking logarithmic derivative of both sides of (3.8)
with respect to $s$ we get
$${\alpha\ln p \over p^s-\alpha }
=h+{\epsilon_p \over s}+\sum_\rho {1\over s-\rho}.
\tag 3.9$$
By letting $s\to\infty$ in (3.9) we find that $h=0$.
Then the stated identity follows from (3.8).

  This completes the proof of the lemma.
\qed\enddemo

  \proclaim{Lemma 3.4}  Let $\alpha, p$ be given in
  Lemma 3.3, and let $s=(1-z)^{-1}$.  Then we have
  $$ {d\over dz}\ln (1-\alpha p^{-s})
  =\sum_{n=0}^\infty \left(\sum_\rho
  [1-(1-1/\rho)^{-n-1}]\right)z^n$$
for $z$ in a small neighborhood of the origin,
where the sum on $\rho$ is over all zeros
of $1-\alpha p^{-s}$
taken in the order given by $|\rho|<T$ for $T\to\infty$.
\endproclaim

 \demo{Proof} By Lemma 3.3 we have
$${d\over dz}\ln (1-\alpha p^{-s})={\epsilon_p\over 1-z}
+\sum_\rho {1\over (1-\rho)+\rho z}{1\over 1-z} \tag 3.10$$
where the sum on $\rho$ is over all nonzero zeros of
$1-\alpha p^{-s}$.  Since
$$\aligned {1\over (1-\rho)+\rho z}{1\over 1-z}
&={1\over 1-\rho}\{\sum_{k=0}^\infty \left({-\rho\over
1-\rho}\right)^kz^k\}\{\sum_{l=0}^\infty z^l\}\\
&=\sum_{n=0}^\infty [1-(1-1/\rho)^{-n-1}] z^n
\endaligned $$
for $z$ in a small neighborhood of the origin,
the stated identity follows from (3.10).
Note that, if $s=\rho=0$ is
a zero of $1-\alpha p^{-s}$, then $(1-1/\rho)^{-n}$ is
interpreted to be $0$ for all positive integers $n$.

  This completes the proof of the lemma.
\qed\enddemo

Remark.  By (3.6), Lemma 3.4, (2.3) and (3.1) we have
$${d\over dz} \log \left(A^{s/2}\xi_N(s)\right)
=\sum_{n=0}^\infty \tau_N(n+1) z^n  $$
with $s=(1-z)^{-1}$, where the $\tau_N(n)$'s are given
in (1.6).   Then, by Lemma 3.2, to prove
that all zeros of $\xi_N(s)$ in the strip $0<\Re s<1$
lie on the critical line $\Re s=1/2$ it is enough to
find an upper bound for each $\tau_N(n)$ which implies
that the above series is analytic for $|z|<1$.

  \proclaim{Lemma 3.5}  Let $\alpha, p$ be given in
  Lemma 3.3, and let $s=(1-z)^{-1}$.  Then we have
  $$ {d\over dz}\ln (1-\alpha p^{-s})
  =\sum_{n=0}^\infty \left(\sum_{j=0}^n \binom {n+1}{j+1}
  {(-1)^j\over j!}\sum_{m=1}^\infty {\ln p\over p^m}\alpha^m
  (\ln p^m)^j \right)z^n$$
for $z$ in a small neighborhood of the origin.
\endproclaim

 \demo{Proof}  Let $m$ be any positive integer.  By using
 mathematical induction on $n$, we can show that
 $${d^n\over dz^n}\left[(1-z)^{-m}p^{-ks}\right]_{|z=0}
 =p^{-k}\sum_{j=0}^n\binom nj (n+m-1)\cdots (j+m)
 (-\ln p^k)^j \tag 3.11 $$
 for $n=1,2,\cdots$.  By using the formula (3.11) with
 $m=2$ we find that
 $$\aligned {d\over dz}\ln (1-\alpha p^{-s})
  &=(1-z)^{-2}\ln p\sum_{k=1}^\infty \alpha^k p^{-ks}\\
  &=\sum_{n=0}^\infty {z^n\over n!}\left(\sum_{k=1}^\infty
 \alpha^k {\ln p\over p^k}\sum_{j=0}^n\binom nj
 (n+1)\cdots (j+2)(-\ln p^k)^j \right)\\
 &=\sum_{n=0}^\infty z^n\left(\sum_{j=0}^n \binom {n+1}{j+1}
  {(-1)^j\over j!}\sum_{k=1}^\infty {\ln p\over p^k}\alpha^k
  (\ln p^k)^j \right).
 \endaligned $$

  This completes the proof of the lemma.
\qed\enddemo

\proclaim{Lemma 3.6}  Let $p$ be a prime, and let
$\alpha$ be a complex number.  Then we have
$$\sum_\rho [1-(1-1/\rho)^{-n}]
  =\sum_{j=1}^n \binom nj {(-1)^{j-1}\over (j-1)!}
\sum_{m=1}^\infty {\ln p\over p^m}\alpha^m(\ln p^m)^{j-1}
  $$
for $n=1,2,\cdots$, where the sum on $\rho$ is over
all zeros of $1-\alpha p^{-s}$ taken in the order given by
$|\Im \rho|<T$ for $T\to\infty$ with a zero of multiplicity
$\ell$ appearing $\ell$ times in the list.
\endproclaim

\demo{Proof}  The stated identity follows from
Lemma 3.4 and Lemma 3.5. \qed\enddemo

  \proclaim{Lemma 3.7}  For $n=1,2,\cdots$, we have
$$\tau_N(n)=\sum_{j=1}^g\tau_{g_j}(n)
+\sum_{l=1}^n \binom nl {(-1)^{l-1}\over (l-1)!}
\sum_{(m, N)>1, m=1}^\infty {\Lambda(m)\over m^{k+1\over 2}}
(\sum_{j=1}^g b_{g_j}(m))(\ln m)^{l-1} $$
where $b_{g_j}(m)$ is given as in (2.2).
\endproclaim

  \demo{Proof}  Let $\tau_H(n)$ be given in (3.3).
By (2.1), (3.1) and (3.6) we have
  $$\tau_N(n)=\tau_H(n)+\sum_{j=1}^g \{
  \sum_{p|N_j}\sum_{\rho_j}[1-(1-1/\rho_j)^{-n}]
  +\sum_{p\nmid N_j, p|N}\sum_{\alpha_j}
  [1-(1-1/\alpha_j)^{-n}]\},\tag 3.12 $$
  where the sum on $\rho_j$ is over all zeros of
$1-\lambda_{g_j}(p)p^{{1-k\over 2}-s}$ with $p|N_j$
and where the sum on $\alpha_j$ is over all zeros of
$1-\lambda_{g_j}(p)p^{{1-k\over 2}-s}+\chi_j(p)p^{-2s}$
 with $p\nmid N_j, p|N$.  By Lemma 3.6 we have
 $$\sum_{\rho_j}[1-(1-1/\rho_j)^{-n}]
 =\sum_{l=1}^n \binom nl {(-1)^{l-1}\over (l-1)!}
\sum_{m=1}^\infty {\ln p\over p^{{k+1\over 2}m}}
b_{g_j}(p^m)(\ln p^m)^{l-1} \tag 3.13 $$
and
$$\sum_{\alpha_j}[1-(1-1/\alpha_j)^{-n}]
=\sum_{l=1}^n \binom nl {(-1)^{l-1}\over (l-1)!}
\sum_{m=1}^\infty {\ln p\over p^{{k+1\over 2}m}}
b_{g_j}(p^m)(\ln p^m)^{l-1}. \tag 3.14 $$
It follows from (3.13) and (3.14) that
$$\aligned &\sum_{p|N_j}\sum_{\rho_j}
  [1-(1-1/\rho_j)^{-n}]+\sum_{p\nmid N_j, p|N}\sum_{\alpha_j}
  [1-(1-1/\alpha_j)^{-n}]\\
  =&\sum_{p|N}\sum_{l=1}^n \binom nl {(-1)^{l-1}\over (l-1)!}
\sum_{m=1}^\infty {\ln p\over p^{{k+1\over 2}m}}
b_{g_j}(p^m)(\ln p^m)^{l-1}. \endaligned \tag 3.15 $$
By (3.15) we have
$$\aligned &\sum_{j=1}^g\left(\sum_{p|N_j}\sum_{\rho_j}
  [1-(1-1/\rho_j)^{-n}]+\sum_{p\nmid N_j, p|N}\sum_{\alpha_j}
  [1-(1-1/\alpha_j)^{-n}]\right)\\
  &=\sum_{l=1}^n \binom nl {(-1)^{l-1}\over (l-1)!}
\sum_{(m, N)>1, m=1}^\infty {\Lambda(m)\over m^{k+1\over 2}}
(\sum_{j=1}^g b_{g_j}(m))(\ln m)^{l-1}.\endaligned $$
The stated identity then follows from (3.12).

  This completes the proof of the lemma.
\qed\enddemo

\heading
4.  Proof of Theorem 1.4
\endheading

   We define an operator $S$ acting on the space $S_k(N, \chi)$ by
  $$\aligned &S(1)=2I \\
  & S(p)=T(p) \\
  &S(p^m)=T(p^m)-\chi(p)p^{k-1}T(p^{m-2})\endaligned  \tag 4.1$$
for $m=2,3,\cdots$.

  \proclaim{Lemma 4.1}  For each prime $p\nmid N$, the trace
$\text{tr}(S(p^m))$ of $S(p^m)$ acting on the space $S_k(N, \chi)$
is given by
$$\text{tr}(S(p^m))=\sum_{j=1}^g b_{g_j}(p^m) \tag 4.2 $$
for $m=0,1,2,\cdots$.
\endproclaim

\demo{Proof}  Let $p$ be any prime.
 It follows from the recursion formula (see (6.25) in Iwaniec [6])
$$T(p^{m+1})=T(p)T(p^m)-\chi(p)p^{k-1}T(p^{m-1})$$
that
$$S(p^{m+1})=S(p)S(p^m)-\chi(p)p^{k-1}S(p^{m-1}) \tag 4.3 $$
for $m=1,2,\cdots$.

  Let $p$ is any prime not dividing $N$.  When $m=0$, we have
$$S(p^m)f_j=2f_j=b_{g_j}(p^m)f_j.$$
When $m=1$, we have
$$S(p^m)f_j=T(p)f_j=\lambda_{g_j}(p)f_j=b_{g_j}(p^m)f_j$$
by Lemma 3.1.  Assume that
$$S(p^m)f_j=b_{g_j}(p^m)f_j$$
for all integers $m\leq l$.  When $m=l+1$, we have
$$\aligned S(p^m)f_j&=(S(p)S(p^l)-\chi(p)p^{k-1}S(p^{l-1}))f_j\\
&=(b_{g_j}(p)b_{g_j}(p^l)-\chi(p)p^{k-1}b_{g_j}(p^{l-1}))f_j\\
&=b_{g_j}(p^m)f_j.\endaligned $$
By mathematical induction the identity
$$S(p^m)f_j=b_{g_j}(p^m)f_j$$
holds for all nonnegative integers $m$.
Since $\{f_1, \cdots, f_g\}$ is a basis for $S_k(N, \chi)$, we have
$$tr(S(p^m))=\sum_{j=1}^g b_{g_j}(p^m)$$
for $m=0, 1,2,\cdots$.

  This completes the proof of the lemma.
\qed\enddemo

  From Lemma 4.1 and the definition (4.1) we
obtain the following corollary.

   \proclaim{Corollary 4.2}  For each prime $p\nmid N$, we have
$$\sum_{j=1}^g b_{g_j}(p^m)=\text{tr}(T(p^m))
-\chi(p)p^{k-1}\,\text{tr}(T(p^{m-2})$$
for $m=0,1,2,\cdots$.
\endproclaim

 \demo{Proof of Theorem 1.4}  By Lemma 3.7 we have
$$\tau_N(n)=\sum_{j=1}^g\tau_{g_j}(n)
+\sum_{l=1}^n \binom nl {(-1)^{l-1}\over (l-1)!}
\sum_{(m, N)>1, m=1}^\infty {\Lambda(m)\over m^{k+1\over 2}}
(\sum_{j=1}^g b_{g_j}(m))(\ln m)^{l-1} \tag 4.4 $$
for $n=1,2,\cdots$.  Since $g_j$ is a normalized
newform in $S_k(N_j,\chi)$, by Theorem 2.1 we have
$$ \aligned \tau_{g_j}(n)=& n\left(\ln {\sqrt {N_j}\over 2\pi}
-\gamma\right)
-\sum_{l=1}^n \binom nl {(-1)^{l-1}\over (l-1)!}\sum_{m=1}^\infty
{\Lambda(m)\over m^{k+1\over 2}} b_{g_j}(m)(\ln m)^{l-1}\\
&+n\left(-{2\over k+1}+\sum_{l=1}^\infty {k+1\over l(2l+k+1)}\right)
+\sum_{m=2}^n\binom nm (-1)^m\sum_{l=1}^\infty
{1\over (l+{k-1\over 2})^m} \endaligned \tag 4.5  $$
By using (4.4) and (4.5) we obtain that
$$\aligned \tau_N(n)&= {n\over 2}\ln(N_1\cdots N_g)
-\sum_{l=1}^n \binom nl {(-1)^{l-1}\over (l-1)!}
\sum_{\underset{(m, N)=1}\to {m=1}}^\infty
{\Lambda(m)\over m^{k+1\over 2}}
\{\sum_{j=1}^g b_{g_j}(m)\}(\ln m)^{l-1}\\
&-ng[\ln 2\pi +\gamma+{2\over k+1}
-\sum_{l=1}^\infty {k+1\over l(2l+k+1)}]
+g\sum_{m=2}^n\binom nm \sum_{l=1}^\infty
{(-1)^m\over (l+{k-1\over 2})^m}. \endaligned  \tag 4.6$$
By the argument in the proof of Theorem 5 in Atkin and Lehner [1], 
we have
$$\prod_{N_j\neq N, 1\leq j\leq g}N_j
=\prod_{m|N, \frak f|m} m^{\nu_m d(N/m)}.\tag 4.7$$
We also have
$$\prod_{N_j= N, 1\leq j\leq g}N_j=N^{\nu_N}. \tag 4.8$$
By (4.6), (4.7) and (4.8) we have
$$\aligned \tau_N(n)=& {n\over 2}
\ln\left(N^{\nu_N}\prod_{\frak f|m, m|N} m^{\nu_md(N/m)}\right)
-ng(\ln 2\pi +\gamma+{2\over k+1})\\
&-\sum_{l=1}^n \binom nl {(-1)^{l-1}\over (l-1)!}
\sum_{\underset{(m, N)=1}\to {m=1}}^\infty
{\Lambda(m)\over m^{k+1\over 2}}
\{\sum_{j=1}^g b_{g_j}(m)\}(\ln m)^{l-1}\\
&+ng\sum_{l=1}^\infty {k+1\over l(2l+k+1)}
+g\sum_{m=2}^n\binom nm \sum_{l=1}^\infty
{(-1)^m\over (l+{k-1\over 2})^m}. \endaligned $$
The stated identity then follows from Corollary 4.2.

This completes the proof of the theorem.
\qed\enddemo

Remark.  If $\chi$ is a primitive Dirichlet character of
modulus $N$, we define
$$L_N(s)=\prod_p\det|1-T(p)p^{-s}+\chi(p)p^{k-1-2s}I|^{-1} $$
where $I$ is the identity map acting on the space $S_k(N, \chi)$
and where the product on $p$ runs over all prime numbers.
Let
$$\xi_N(s)=N^{gs/2}(2\pi)^{-gs}\Gamma^g({k-1\over 2}+s)
L_N({k-1\over 2}+s), $$
where $g$ denotes the dimension of the space $S_k(N,\chi)$.
Since the space $S_k(N,\chi)$ contains only newforms (see
\S6.7 of Iwaniec [6]), we have $\xi_N(s)=\xi_H(s)$ by
(3.6).  Hence, $\xi_N(s)$ is an entire function and satisfies
the functional identity
$$\xi_N(s)=w\bar\xi_N(1-\bar s),$$
where $w$ is a constant.

\enddocument